\documentclass[12pt,thmsa]{article}
\usepackage{amssymb,amsmath}

\usepackage{sw20bams}

\numberwithin{equation}{section}

\newcounter{saveeqn}
\newcommand{\alphaeqn}{\setcounter{saveeqn}{\value{equation}}%
\setcounter{equation}{0}%
\global\def\theequation{\mbox{\thesection.\arabic{saveeqn}\alph{equation}}}}
\newcommand{\reseteqn}{\setcounter{equation}{\value{saveeqn}}%
\global\def\theequation{\thesection.\arabic{equation}}}

\newcommand{\QED}{\nobreak\strut\hfill$\square$\par\vskip3mm}

\def\jj{\ell}
\def\jjj{\ell}
\def\iiii{r}
\def\jjjj{s}
\def\ttt{\ell}
\def\tttt{\beta}
\def\ssss{\alpha}
\def\rrr{r}
\def\rrrr{t}
\def\ttttt{s}
\def\sssss{r}

\begin{document}

\title{SOME NEW FORMULAS FOR $\pi $}
\author{Gert Almkvist, Christian 
Krattenthaler$^{\dagger}$\thanks{Research partially supported by the Austrian
Science Foundation FWF, grant P13190-MAT, and by EC's IHRP Programme,
grant HPRN-CT-2001-00272.}
\hbox{\ }and Joakim
Petersson\\[12pt]
\rm
Matematiska Institutionen, Lunds Universitet, \\
\rm
SE-221 00 Lund, Sweden.\\
\rm E-mail: {\tt Gert.Almkvist@math.lu.se}, 
{\tt Joakim.Petersson@math.lu.se}\\[6pt]
\rm
$^\dagger$Institut f\"ur Mathematik der Universit\"at Wien,\\
\rm Strudlhofgasse 4, A-1090 Wien, Austria.\\
\rm E-mail: {\tt KRATT@Ap.Univie.Ac.At}\\
\rm WWW: \tt http://www.mat.univie.ac.at/\~{}kratt}
\date{}

\maketitle

\begin{abstract}
We show how to find 
series expansions for $\pi$ of the
form $\pi=\sum_{n=0}^\infty {S(n)}\big/{\binom{mn}{pn}a^n}$, where
$S(n)$ is some polynomial in $n$ (depending on $m,p,a$). 
We prove that there exist such
expansions for $m=8k$, $p=4k$, $a=(-4)^k$, for any $k$, and give
explicit examples for such expansions for small values of $m$, $p$ and $a$.
\end{abstract}

\bigskip
\section{Introduction}
\bigskip

Using the formula (due to Bill Gosper \cite{Gosper})
\[
\pi =\sum_{n=0}^\infty \frac{50n-6}{\binom{3n}n2^n} ,
\]
Fabrice Bellard \cite[file {\tt pi1.c}]{Bellard} 
found an algorithm for computing the $n$-th decimal of $\pi $ without
computing the earlier ones. Thus he improved an earlier algorithm due to
Simon Plouffe \cite{Plouffe}. 

This formula can be proved in the following way.
We have
\[
\frac 1{\binom{3n}n}=(3n+1)\int_0^1x^{2n}(1-x)^ndx .
\]
Hence the right hand side of the formula will be
\[
\int_0^1\sum_{n=0}^\infty (50n-6)(3n+1)\left(\frac {x^2(1-x)}{2}\right)^ndx
.
\]
But
\[
\sum_{n=0}^\infty (50n-6)(3n+1)y^n=\frac{2(56y^2+97y-3)}{(1-y)^3} ,
\]
so we get
\begin{align*}
RHS&=8\int_0^1\frac{28x^6-56x^5+28x^4-97x^3+97x^2-6}{(x^3-x^2+2)^3}dx\\
&=
\left[ \frac{4x(x-1)(x^3-28x^2+9x+8)}{(x^3-x^2+2)^2}+4\arctan (x-1)\right]
_0^1=\pi .
\end{align*}
We then asked if there are other such formulas of type
\[
\pi =\sum_{n=0}^\infty \frac{S(n)}{\binom{mn}{pn}a^n} ,
\]
where $S(n)$ is a polynomial in $n$. Using the same trick as before we have
\[
\frac 1{\binom{mn}{pn}}=(mn+1)\int_0^1x^{pn}(1-x)^{(m-p)n}\,dx ,
\]
and hence
\begin{equation} \label{eq:P1}
RHS=\int_0^1\sum_{n=0}^\infty (mn+1)S(n)\left(\frac 
{x^p(1-x)^{m-p}}{a}\right)^n\,dx .
\end{equation}
If $S(n)$ has degree $d$ then
\begin{equation} \label{eq:T(y)} 
\sum_{n=0}^\infty (mn+1)S(n)y^n=\frac{T(y)}{(1-y)^{d+2}} ,
\end{equation}
where $T$ has degree $d+1$. It follows with
\[
y=\frac{x^p(1-x)^{m-p}}a 
\]
that
\begin{equation} \label{eq:P2}
RHS=\int_0^1\frac{P(x)}{(x^p(1-x)^{m-p}-a)^{d+2}}\,dx ,
\end{equation}
where $P(x)$ is a polynomial in $x$ of degree $m(d+1)$. 
We want this integral to be equal to $\pi $.
A good way to get $\pi $ is to have $\arctan (x)$ or arctan$(x-1)$ after
integration. This means that
\[
x^p(1-x)^{m-p}-a 
\]
must have the factor $x^2+1$ or $(x-1)^2+1$, that is have a zero at $i$ or $%
1+i$. This restricts $m$ and $p$ and gives the value of $a$. After
experimenting with the LLL-algorithm we found formulas for $\pi $ in the
following cases (there could be many more):

\smallskip

\bigskip

\bigskip

\begin{center}
\begin{tabular}{c|c|l|cl}
$m$ & $p$ & $\hphantom{-}a$ & $\deg(S)$ &  \\ 
\cline{1-4}
\hphantom{1}3 & \hphantom{1}1 & $\hphantom{-}2$ & \hphantom{1}1 & Gosper 
\cite{Gosper}\\ 
\hphantom{1}7 & \hphantom{1}2 & $\hphantom{-}2$ & \hphantom{1}5 & Bellard 
\cite{Bellard}$^1$\\ 
\hphantom{1}8 & \hphantom{1}4 & $-4$ & \hphantom{1}4 &  \\ 
10 & \hphantom{1}4 & $\hphantom{-}4$ & \hphantom{1}8 &  \\ 
12 & \hphantom{1}4 & $-4$ & \hphantom{1}8 &  \\ 
16 & \hphantom{1}8 & $\hphantom{-}16$ & \hphantom{1}8 &  \\ 
24 & 12 & $-64$ & 12 &  \\ 
32 & 16 & $\hphantom{-}256$ & 16 &  \\ 
40 & 20 & $-4^5$ & 20 &  \\ 
48 & 24 & $\hphantom{-}4^6$ & 24 &  \\ 
56 & 28 & $-4^7$ & 28 &  \\ 
64 & 32 & $\hphantom{-}4^8$ & 32 &  \\ 
72 & 36 & $-4^9$ & 36 &  \\ 
80 & 40 & $\hphantom{-}4^{10}$ & 40 & 
\end{tabular}
\end{center}

\def\thefootnote{}
\footnote{$^1$ This formula was used in Bellard's world
record setting
computation of the 1000 billionth {\it binary} digit of $\pi$, being based
on the algorithm in \cite{BaBoPl}.}

\bigskip

E.g., to find the formula in the last case we computed
\[
s(k)=\sum_{n=0}^\infty \frac{n^k}{\binom{80n}{40n}4^{10n}} 
\]
for $k=0,1,2,\dots, 40$ with $6000$ digits. Then the
LLL-algorithm found (in about three days) the linear combination between $%
\pi $ and $s(0),s(1),\dots ,s(40)$. A good check is that the
coefficient of $\pi $ is a product of small primes (actually all less than $%
80$).

\bigskip
\section{Proving the formulas}
\bigskip

For $m\leqslant 16$ the integrals can be computed by brute force using
{\sl Maple}. 
The higher cases are all symmetric, i.e., $m=2p$. Using the symmetry
\[
x\longleftrightarrow 1-x 
\]
one can assume that $\arctan (x)$ and $\arctan (x-1)$ each contribute 
$\pi /2$. So we make the wild assumption that
\[
\int \frac{P(x)}{Q(x)^{d+2}}\,dx=\frac{R(x)}{Q(x)^{d+1}}+2\arctan (x)+2\arctan
(x-1) ,
\]
where $R(x)$ is a polynomial with $R(0)=R(1)=0$. Differentiation with
respect to $x$ yields
\[
\frac P{Q^{d+2}}=\frac{R^{\prime }}{Q^{d+1}}-(d+1)\frac{Q^{\prime }R}{Q^{d+2}%
}+2\left(\frac 1{x^2+1}+\frac 1{x^2-2x+2}\right), 
\]
or
\[
QR^{\prime }-(d+1)Q^{\prime }R=P-2Q^{d+2}\left(\frac 1{x^2+1}+
\frac 1{x^2-2x+2}\right) .
\]
This is a differential equation of first order in $R$ with polynomial
coefficients. Taking the case $\binom{32n}{16n}$, 
we have $\deg (P)=544$ and $d=16$, 
so the denominator $Q^{d+2}$ has degree $576$. We solve the
differential equation using {\sl Maple}'s
\[
\text{\tt dsolve(*,R(0)=0,series}) 
\]
setting the order to $600$. Then we find that $R$ is a polynomial of degree 
$543$. In practice one replaces the coefficients of $\frac 1{x^2+1}$ and 
$\frac 1{x^2-2x+2}$ by $2$ times the denominator in the formula found for 
$\pi $. Then $R$ will have (huge) integer coefficients.

\begin{remark}\rm
Computing the integral in the case $\binom{8n}{4n}$, 
{\sl Maple} got the terms
\[
\arctan\left( \frac {2x^3-3x^2+7x-3}{5}\right)+\arctan 
\left(\frac x2-\frac 14\right) ,
\]
which can be shown to be equal to
\[
\arctan (x)+\arctan (x-1).
\]
In the nonsymmetric case $\binom{12n}{4n}$ one gets $\arctan (x-1)$ and $%
\arctan (x^3-2x^2+x-1)$, but the latter term makes no contribution since it
has the same values at $0$ and $1$. 
\end{remark}

\begin{remark}\rm
When we set this up then,
in principle, we could let {\sl Maple} do the evaluation of $T(y)$
(and thus of $P(x)$) through \eqref{eq:T(y)}, 
i.e., let {\sl Maple}'s formal summation
routines compute the sum
\[
\sum_{n=0}^\infty  (mn+1)S(n)y^n,
\]
or, to make the point clearly visible, compute a sum of the form
$$
\sum_{n=0}^\infty \tilde S(n) y^n,
$$
where $\tilde S(n)$ is a certain polynomial in $n$. However, it turns out
that this is very time-consuming. It is a much better idea to expand 
$\tilde S(n)$ as 
$\tilde S(n)=\sum _{j=0} ^{d+1}a(j)\binom {n+j}{d+1}$, and then
use that $\sum _{n=0} ^{\infty}\binom {n+j}{d+1}y^n=y^{d+1-j}/(1-y)^{d+2}$
for $j=0,1,\dots,d+1$. See also the approach that we follow in
Section~\ref{sec:Proof} when we prove our main result. 
\end{remark}

\bigskip
\section{The case $\dbinom{8kn}{4kn}$, $k=1,2,3,\dots $}
\bigskip

Solving the differential equation for $R(x)$ one finds that $R(x)$ has the
factor\break $x(1-x)(2x-1)$. If one divides out the factor $2x-1$, 
the remaining
factor is invariant under the substitution $x\longleftrightarrow 1-x$.
Hence we make the substitution
\[
R(x)=(2x-1)\check{R}(t) ,
\]
where
\[
z=x(1-x) .
\]
Then
\[
\frac{dz}{dx}=1-2x\text{ and }(1-2x)^2=1-4z .
\]
Now
\[
Q=z^{4k}-(-4)^k ,
\]
so
\[
Q^{\prime }=4kz^{4k-1}(1-2x) 
\]
and
\[
R^{\prime }=2\check{R}(z)+(2x-1)\frac{d\check{R}}{dz}(1-2x)=2\check{R}-(1-4z)%
\frac{d\check{R}}{dz} .
\]
It follows that the left hand side of the differential equation is
\begin{align*} 
QR^{\prime }-(4k+1)Q^{\prime }R&=Q\left(2\check{R}-(1-4z)\frac{d\check{R}}{dz}%
\right)+(4k+1)4kz^{4k-1}(1-2x)^2\check{R}\\
&=(2Q+4k(4k+1)(1-4z)z^{4k-1})\check{R}-(1-4z)Q\frac{d\check{R}}{dz}. 
\end{align*}
On the right hand side we have
\[
\frac 1{x^2+1}+\frac 1{x^2-2x+2}=\frac{3-2z}{z^2-2z+2} ,
\]
and $Q$ is for all $k$ divisible by
\[
z^4+4=(z^2-2z+2)(z^2+2z+2) .
\]
Hence $P$ and $Q$ and
\[
Q\left(\frac 1{x^2+1}+\frac 1{x^2-2x+2}\right) 
\]
are all polynomials in $z$. Finally we get the following
differential equation for $\check{R}$:
\begin{equation} \label{eq:diff}
-(1-4z)Q\frac{d\check{R}}{dz}+(2Q+4k(4k+1)(1-4z)z^{4k-1})\check{R}-P+2(3-2z)%
\frac{Q^{4k+2}}{z^2-2z+2}=0 .
\end{equation}

If $P$ is known then one solves this equation for $\check{R}$ just as before
(but the degrees are cut in half). 

But there is a possibility to {\it
find\/} the
formula for $\pi $ {\it and prove} it {\it in one stroke}. Let
\[
N(k)=4k(4k+1) ,
\]
and assume that
\begin{equation} \label{eq:Ansatz1}
\check{R}(z)=\sum_{j=1}^{N(k)-1}a(j)z^j 
\end{equation}
and
\begin{equation} \label{eq:pi}
\pi =\sum_{n=0}^\infty \frac{S_k(n)}{\binom{8kn}{4kn}(-4)^{kn}} ,
\end{equation}
with
\begin{equation} \label{eq:Ansatz2}
S_k(n)=\sum_{j=0}^{4k}a(N(k)+j)n^j .
\end{equation}
Substituting this into the differential equation we will get a system of
linear equations for the $a(j)$'s of size $(N(k)+4k)\times (N(k)+4k)$. 
If this system is nonsingular, we can be sure to be able to solve it
and thus {\it find and prove} an expansion for $\pi$ of the form
\eqref{eq:pi}.

We went (again) to the computer and generated the system of equations
for small values of $k$.
Let $A(k)$ be the matrix of the system. 
Let further $r(k)$ be the least common
denominator of the coefficients of $S_k(n)$. 
We obtained the following tables:

\bigskip

\begin{tabular}{l|l}
\hphantom{1}$k$ & $\hphantom{-}\det(A(k))$ \\ 
\hline\\[-8pt]
\hphantom{1}1 & $\hphantom{-}2^{91}3^85^77^2$ \\ 
\hphantom{1}2 & $-2^{523}3^{52}5^{17}7^{14}11^413^3$ \\ 
\hphantom{1}3 & $\hphantom{-}2^{1367}3^{177}5^{41}7^{25}11^{20}13^{19}17^519^423^2$ \\ 
\hphantom{1}4 & $-2^{3231}3^{167}5^{83}7^{53}11^{28}13^{27}17^{25}19^823^629^331^2$ \\ 
\hphantom{1}5 & $\hphantom{-}
2^{5399}3^{290}5^{345}7^{93}11^{41}13^{37}17^{33}19^{32}23^{10}29^731^637^3$%
\end{tabular}
\bigskip

\begin{tabular}{r|l}
$k$ & $\hphantom{-}r(k)$ \\ 
\hline\\[-8pt]
1 & $\hphantom{-}3^25^27^2$ \\ 
2 & $\hphantom{-}3^65^37^211^213^2$ \\ 
3 & $\hphantom{-}2^53^35^27^211^213^217^219^223^2$ \\ 
4 & $\hphantom{-}2^33^{10}5^67^311\cdot 13^217^219^223^229^231^2$ \\ 
5 & $\hphantom{-}2^{11}3^{13}5^37^411^313^317^219^223^229^231^237^2$ \\ 
6 & $\hphantom{-}2^{10}3^75^87^311^313^317^219^223^229^231^237^241^243^247^2$ \\ 
7 & $\hphantom{-}2^{17}3^{16}5^97^311^413^317^319^223^229^231^237^241^243^247^253^2$ \\ 
8 & $\hphantom{-}
2^{10}3^{19}5^87^811^413^317^319^323^229^231^237^241^243^247^253^259^261^2$
\\ 
9 & $\hphantom{-}
2^{23}3^45^{10}7^811^413^417^319^323^329^231^237^241^243^247^253^259^261^267^271^2 
$ \\ 
10 & $\hphantom{-}
2^{20}3^{20}5^47^911^513^417^319^323^229^231^237^241^243^247^253^259^261^267^271^273^279^2 
$%
\end{tabular}

\bigskip

{}From the table it seems ``obvious" that the determinant of $A(k)$
will never vanish as it grows rather quickly in absolute value with $k$. 
We also see that the same primes occur in the factorizations of $r(k)$ and $\det
(A(k))$. Even more striking, it seems that the largest prime factors
occurring grow only slowly (namely approximately linearly) 
when $k$ increases. The last fact strongly indicates that there may
even be a closed form formula for $\det(A(k))$. As it turns out, this
is indeed the case. In the next section we will explicitly compute
the determinant of an equivalent system of equations, from which it
follows that (see the remark after Theorem~\ref{T2})
\begin{equation} \label{eq:det(A(k))}
\det (A(k))=(-1)^{k-1}
2^{32k^3+24k^2+2k-1}k^{8k^2+2k}((4k+1)!)^{4k}\frac{(8k)!}{(4k)!}
\prod_{j=1}^{4k}\frac{(2j)!}{j!}.
\end{equation}
Hence we have the following theorem, which is the main theorem of our
paper.

\begin{theorem} \label{T1}
For all $k\geq 1$ there is a formula
\[
\pi =\sum_{n=0}^\infty \frac{S_k(n)}{\binom{8kn}{4kn}(-4)^{kn}},
\]
where $S_k(n)$ is a polynomial in $n$ 
of degree $4k$ with rational coefficients. The polynomial $S_k(n)$ can be
found by solving the system of linear equations generated by
\eqref{eq:diff} and the Ansatz \eqref{eq:Ansatz1} and
\eqref{eq:Ansatz2}.
\end{theorem}

The denominators of all $a(j)$'s divide $r(k)$ which is much smaller than 
$\det (A(k))$. This means that there must be some miracle occurring at the end
when solving the system. E.g., when $k=5$ then $\det (A(5))$ has about $2400$
digits but $r(5)$ only $40$. Unfortunately, we are not able to
offer an explanation for that.

In practice we are only interested in the
coefficients of $S_k$ so we try to eliminate $a(1),a(2),\dots
,a(N(k)-1)$ first. This can be done by first avoiding all equations
containing $a(N(k)),\dots ,a(N(k)+4k)$, i.e., the equations
coming from the coefficients of $t^{4kv}$, $v=0,1,\dots ,4k+1$.
This looks very nice theoretically, in particular as the system for $%
a(1),\dots ,a(N(k))$ is triangular, but in practice the computer
breaks down since the rational numbers occuring become very large. 

\medskip
We close this section by listing a few explicit examples.

\begin{example}
We have
\[
\pi =\frac 1r\sum_{n=0}^\infty \frac{S(n)}{\binom{8n}{4n}(-4)^n} ,
\]
where
\[
r=3^25^27^2 
\]
and
\[
S(n)=-89286+3875948n-34970134n^2+110202472n^3-115193600n^4 .
\]
\end{example}

\begin{example}
We have
\[
\pi =\frac 1r\sum_{n=0}^\infty \frac{S(n)}{\binom{16n}{8n}16^n}, 
\]
where
\[
r=3^65^37^211^213^2 
\]
and
\begin{multline*} 
S(n)=-869897157255-3524219363487888n+112466777263118189n^2 
\\
-1242789726208374386n^3+6693196178751930680n^4-19768094496651298112n^5 
\\
+32808347163463348736n^6-28892659596072587264n^7+10530503748472012800n^8. 
\end{multline*}
\end{example}

\begin{example}
We have
\[
\pi =\frac 1r\sum_{n=0}^\infty \frac{S(n)}{\binom{32n}{16n}256^n} ,
\]
where
\[
r=2^33^{10}5^67^311\cdot 13^217^219^223^229^231^2 
\]
and
\begin{align*}
S(n)=&-2062111884756347479085709280875 
\\
&+1505491740302839023753569717261882091900n 
\\
&-112401149404087658213839386716211975291975n^2 
\\
&+3257881651942682891818557726225840674110002n^3 
\\
&-51677309510890630500607898599463036267961280n^4 
\\
&+517337977987354819322786909541179043148522720n^5 
\\
&-3526396494329560718758086392841258152390245120n^6 
\\
&+171145766235995166227501216110074805943799363584n^7 
\\
&-60739416613228219940886539658145904402068029440n^8 
\\
&+159935882563435860391195903248596461569183580160n^9 
\\
&-313951952615028230229958218839819183812205608960n^{10} 
\\
&+457341091673257198565533286493831205566468325376n^{11} 
\\
&-486846784774707448105420279985074159657397780480n^{12} 
\\
&+367314505118245777241612044490633887668208926720n^{13} 
\\
&-185647326591648164598342857319777582801297080320n^{14} 
\\
&+56224688035707015687999128994324690418467340288n^{15} 
\\
&-7687255778816557786073977795149360408612044800n^{16} .
\end{align*}

\end{example}
\bigskip

\begin{example}
(nonsymmetric). We have
\[
\pi =\frac 1r\sum_{n=0}^\infty \frac{S(n)}{\binom{10n}{4n}4^n} ,
\]
where 
\[
r=3^25\cdot 7^211\cdot 13\cdot 17\cdot 19 
\]
and
\begin{multline*}
S(n)=-4843934523072-1008341177146848n+23756198610834352n^2\\
-242873913552020704n^3 
+1195813551184400032n^4-3272960363556054592n^5\\
+4909379167837011328n^6 
-3816399750842818816n^7+1190182007407360000n^8 .
\end{multline*}
\end{example}

\bigskip
\section{Proof of the Theorem}
\label{sec:Proof}
\bigskip

We want to prove that, by making the Ansatz \eqref{eq:Ansatz1} and
\eqref{eq:Ansatz2} and substituting this into the differential equation
\eqref{eq:diff} (the polynomial $P$ being given by $S_k$ through
\eqref{eq:P1}--\eqref{eq:P2}, $Q$ being given by
$(x(1-x))^{4k}-(-4)^{k}$), 
the resulting system of linear equations will always have a solution. 
In fact, we aim at finding an explicit formula for the determinant of
the corresponding matrix of coefficients that allows us to conclude
that it can never vanish.

It turns out that for that purpose it is more convenient 
to set up the system of linear equations in a different, but
equivalent way. This equivalent system will have coefficient matrix
$M$ (see \eqref{eq:defM}). The evaluation of its determinant will
be accomplished through Eqs.~\eqref{eq:M'}, 
\eqref{eq:M''}, \eqref{eq:M'''}, \eqref{eq:M^X}, and Theorem~\ref{T2}.

To be precise, we encode the polynomial $S_k(x)$ (and, thus,
$T(y)$ and $P(x)$) differently.
We claim that $T(y)$ has an expansion of the form 
\begin{equation} \label{eq:A3}
T(y)=\sum _{j=0}
^{4k+1}c_jy^j,
\end{equation} 
subject to the single constraint
\begin{equation} \label{eq:CK1}
\prod _{i=1} ^{4k+1}(4ik-1)c_0+\sum _{j=1} ^{4k+1}(-1)^j
\left(\prod _{i=1} ^{4k+1-j}(4ik-1)\right)
\left(\prod _{i=1} ^{j-1}(4ik+1)\right)c_j\,=\,0.
\end{equation}
(As usual, empty poducts have to be interpreted as 1.)
This is seen as follows. The polynomial $S_k(x)$ can be written in the
form 
$$S_k(n)=\sum _{j=0} ^{4k}\binom {n+j}{4k}s_j,$$
for some coefficients $s_j$.
Hence, we have
\begin{align*} &\frac {T(y)} {(1-y)^{4k+2}}=\sum _{n=0} ^{\infty}
(4kn+1)S_k(n)y^n\\
&\kern15pt=\sum _{n=0} ^{\infty}(4kn+1)\sum _{j=0} ^{4k}\binom
{n+j}{4k}s_jy^n\\
&\kern15pt=\sum _{j=0} ^{4k}s_j\sum _{n=0} ^{\infty}(4k(n+j+1)-4k(j+1)+1)\binom
{n+j}{4k}y^n\\
&\kern15pt=\sum _{j=0} ^{4k}s_j\sum _{n=0} ^{\infty}
\left(4k(4k+1)\binom {n+j+1}{4k+1}y^n-(4k(j+1)-1)\binom {n+j}{4k}y^n\right)\\
&\kern15pt=\sum _{j=0} ^{4k}s_j\left(4k(4k+1)\frac {y^{4k-j}} {(1-y)^{4k+2}}
-(4k(j+1)-1)\frac {y^{4k-j}} {(1-y)^{4k+1}}\right)\\
&\kern15pt=\frac {1} {(1-y)^{4k+2}}\sum _{j=0}
^{4k+1}y^j\left((4jk+1)s_{4k-j}+
(4k(4k-j+2)-1)s_{4k-j+1}\right).
\end{align*}
In the last line, $s_{-1}$ and $s_{4k+1}$ have to be read as $0$.
It is now a trivial exercise to substitute the coefficients of $y^j$
in the sum in the last line into the left-hand side of \eqref{eq:CK1} and verify the truth
of \eqref{eq:CK1}.

The above implies that
$$P(z)=\sum _{j=0} ^{4k+1}c_jz^{4jk}/(-4)^{jk},$$
where the coefficients $c_j$ obey \eqref{eq:CK1}. 

Now we are ready to set up the system of linear equations. We make
again the Ansatz \eqref{eq:Ansatz1}, but we replace \eqref{eq:Ansatz2} by
\begin{equation} \label{eq:P(t)}
P(z)=\sum _{j=0} ^{4k+1}a(N(k)+j)z^{4jk}/(-4)^{jk},
\end{equation}
where the $a(N(k)+j)$, $j=0,1,\dots,4k+1$, are subject to \eqref{eq:CK1} (i.e.,
the relation \eqref{eq:CK1} holds when $c_j$ is replaced by $a(N(k)+j)$).
Clearly, we have to add \eqref{eq:CK1} to the set of equations that result from the
differential equation \eqref{eq:diff}.

The coefficient matrix of the system looks a follows:
\begin{equation} \label{eq:defM}
M=\begin{pmatrix} x&y\\U&V\end{pmatrix},
\end{equation}
where $x$ is a line vector of $N(k)-1$ zeroes,
$y=(y_0,y_1,\dots,y_{4k+1})$ is the vector of
coefficients of \eqref{eq:CK1}, i.e., $y_0=\prod _{i=1} ^{4k+1}(4ik-1)$,
and 
$$y_{\ell}=(-1)^\ell\left(\prod _{i=1} ^{4k+1-\ell}(4ik-1)\right)
\left(\prod _{i=1} ^{\ell-1}(4ik+1)\right),$$ 
$\ell=1,2,\dots,4k+1$,
where $U$ is an $(N(k)+4k)\times (N(k)-1)$ matrix and $V$ is an
$(N(k)+4k)\times (4k+2)$ matrix, both of which we define below.

We consider the top-most line of $M$ (which is formed out of $x$ and
$y$) as row~0 of $M$.
We label the rows of $U$ and $V$ by $i$ running from $1$ to $N(k)+4k$.
Furthermore, we label the columns of $M$ by $j$ running from $1$ to
$N(k)+4k+1$. 

Following this labelling scheme,
the matrix $U$ has nonzero entries only in the four diagonals $i=j$,
$i=j+1$, $i=j+4k$, $i=j+4k+1$. We denote the entries in column~$j$ on these
four diagonals in order $f_0(j)$, $f_1(j)$, $g_0(j)$, and $g_1(j)$,
where
\begin{align*} f_0(j)&=j(-4)^k,\\
f_1(j)&=-(4j+2)(-4)^k,\\
g_0(j)&=(N(k)-j),\\
g_1(j)&=-(4N(k)-4j-2).
\end{align*}
To be precise, the $(j,j)$-entry is $f_0(j)$,
the $(j+1,j)$-entry is $f_1(j)$,
the $(j+4k,j)$-entry is $g_0(j)$,
the $(j+4k+1,j)$-entry is $g_1(j)$, $j=1,2,\dots,N(k)-1$.

On the other hand, the matrix $V$ is composed out of columns, labelled
$N(k),\break N(k)+1,\dots,N(k)+4k+1$, each of which containing just one
nonzero entry. To be precise, the nonzero entry of column $N(k)+j$ is
located in the $(4jk+1)$-st row (according to our labelling scheme),
and it is equal to $(-4)^{-jk}$, $j=0,1,\dots,4k+1$.

\medskip
We will now compute the determinant of $M$ and show that it does not
vanish for any $k$. 

\medskip
We perform some row operations on $M$, with the effect that the
entries of $y$ get eliminated. This is achieved by subtracting 
$\prod _{i=1} ^{4k+1}(4ik-1)$ times row~1 from row~0, and
$$(-1)^\jjj (-4)^{\jjj k}\left(\prod _{i=1} ^{4k+1-\jjj }(4ik-1)\right)
\left(\prod _{i=1} ^{\jjj -1}(4ik+1)\right)$$
times row $4\jjj k+1$ from row~0, $\jjj =1,2,\dots,4k+1$.
Doing this, we must expect changes in row~0 in columns $1$,
$4k$, $4k+1$, $8k$, $8k+1$, \dots,$N(k)-4k=16k^2$. However, at this
point a miracle occurs: the new entries in row~0 in columns $4\jjj k+1$,
$\jjj =0,1,\dots, 4k-1$, are still 0. On the other hand, 
the values of the new entries in
row~0 in columns $4\jjj k$, $\jjj =1,2,\dots,4k$ are
\begin{equation} \label{eq:CK2}
(-1)^{\jjj -1}(-4)^{(\jjj +1)k}8k(4k+1)\left(\prod _{i=1} ^{4k-\jjj }(4ik-1)\right)
\left(\prod _{i=1} ^{\jjj -1}(4ik+1)\right).
\end{equation}

After these manipulations we obtain a matrix of the form
$$M'=\begin{pmatrix} x'&y'\\U&V\end{pmatrix},$$
where $x'$ and $y'$ have the same dimensions as before $x$ and $y$,
respectively, where the only nonzero entries of $x'$ are in columns
labelled by numbers which are divisible by $4k$, with the entry in
column $4\jjj k$ given by \eqref{eq:CK2}, and where $y'$ consists only of zeroes.
We have 
\begin{equation} \label{eq:M'}
\det M=\det M'.
\end{equation}

The next step consists in expanding the determinant of $M'$ with
respect to columns $N(k),N(k),\dots,N(k)+4k+1$ (i.e., the last $4k+2$
columns). Since each of these columns contains just one nonzero entry
(which is a power of $(-4)^{-k}$),
we have
\begin{equation} \label{eq:M''}
\det M'=\pm(-4)^{-k(2k+1)(4k+1)}\det M'',
\end{equation}
where $M''$ is the matrix arising from $M'$ by deleting the last
$4k+2$ columns and the rows $4\jjj k+1$, $\jjj =0,1,\dots,4k+1$.
More precisely, the matrix $M''$ has the following form:
$$\begin{pmatrix} x_1&x_2&x_3&\hdotsfor2 &x_{4k+1}\\
F_1&0&0&\hdotsfor2 &0\\
G_1&F_2&0&\hdotsfor2 &0\\
0&G_2&F_3&\hdotsfor2 &0\\
0&0&G_3&\ddots& &\vdots\\
 & & \ddots&\ddots&&0\\
&& &0 & G_{4k}&F_{4k+1}\\
& & &\dots & 0&G_{4k+1}
\end{pmatrix},
$$
where $x_\jjj $, $\jjj =1,2,\dots,4k$, is a line vector with $4k$ entries, all
of them being zero except for the last, which is equal to \eqref{eq:CK2},
where $x_{4k+1}$ is a line vector of $4k-1$
zeroes, where $F_\jjj $ and $G_\jjj $, $\jjj =1,2,\dots,4k$, are $(4k-1)\times
(4k)$ matrices with nonzero entries only in the two main diagonals,
and where $F_{4k+1}$ and $G_{4k+1}$ are $(4k-1)\times (4k-1)$
matrices, $G_{4k+1}$ being upper triangular.
To be precise, for $\jjj =1,2,\dots,4k$ we have
$$
F_\jjj =\left(\smallmatrix f_1(4(\jjj -1)k+1)&f_0(4(\jjj -1)k+2)&0&\dots\\
0&f_1(4(\jjj -1)k+2)&f_0(4(\jjj -1)k+3)&0&\dots\\
&\ddots&\ddots&\\
&&\ddots&\ddots&\\
&&0&f_1(4\jjj k-2)&f_0(4\jjj k-1)&0\\
&&&0&f_1(4\jjj k-1)&f_0(4\jjj k)\\
\endsmallmatrix\right)
$$
and
$$
G_\jjj =\left(\smallmatrix g_1(4(\jjj -1)k+1)&g_0(4(\jjj -1)k+2)&0&\dots\\
0&g_1(4(\jjj -1)k+2)&g_0(4(\jjj -1)k+3)&0&\dots\\
&\ddots&\ddots&\\
&&\ddots&\ddots&\\
&&0&g_1(4\jjj k-2)&g_0(4\jjj k-1)&0\\
&&&0&g_1(4\jjj k-1)&g_0(4\jjj k)\\
\endsmallmatrix\right),
$$
and we have
$$G_{4k+1}=\left(\smallmatrix g_1(16k^2+1)&g_0(16k^2+2)&0&\dots\\
0&g_1(16k^2+2)&g_0(16k^2+3)&\dots\\
&\ddots&\ddots&\\
&&\ddots&\ddots&\\
&&&g_1(16k^2+4k-2)&g_0(16k^2+4k-1)\\
&&&0&g_1(16k^2+4k-1)\\
\endsmallmatrix\right).
$$
The precise form of $F_{4k+1}$ is without relevance for us.
We do a Laplace expansion with respect to the last $4k-1$ rows.
Because of the triangular form of $G_{4k+1}$ we obtain
\begin{equation} \label{eq:M'''}
\det M''=\left(\prod _{i=16k^2+1} ^{16k^2+4k-1}g_1(i)\right)\det M''',
\end{equation}
where 
$$M'''=\begin{pmatrix} x_1&x_2&x_3&\dots&x_{4k}\\
F_1&0&0&\dots\\
G_1&F_2&0&\dots\\
0&G_2&F_3&\dots\\
0&0&G_3&\ddots&\vdots\\
&\ddots &\ddots&\ddots&0\\
 & & &G_{4k-1}&F_{4k}\\
& & &0 & G_{4k}
\end{pmatrix}.
$$

Instead of $M'''$ we consider a more general matrix. Define the functions
\begin{align*} f_0(t,j)&=((N(k)+j)Y_t-X_{2,t})(-4)^k,\\
f_1(t,j)&=-((4N(k)+4j+2)Y_t-4X_{1,t})(-4)^k,\\
g_0(t,j)&=(X_{2,t}-jY_t),\\
g_1(t,j)&=-(4X_{1,t}-(4j+2)Y_t).
\end{align*}
It should be noted that these functions specialize to
$f_0(j),f_1(j),g_0(j),g_1(j)$, respectively, 
if $X_{1,t}=X_{2,t}=N(k)$ and $Y_t=1$.
Now we define the matrix $M^X$ by
\begin{equation} \label{eq:CK3}
M^X=\begin{pmatrix} x_1&x_2&x_3&\dots&x_{4k}\\
F^X_1&0&0&\dots\\
G^X_1&F^X_2&0&\dots\\
0&G^X_2&F^X_3&\dots\\
0&0&G^X_3&\ddots&\vdots\\
&\ddots &\ddots&\ddots&0\\
 & & &G^X_{4k-1}&F^X_{4k}\\
& & &0 & G^X_{4k}
\end{pmatrix},
\end{equation}
where
$$
F^X_\jjj =\left(\smallmatrix f_1(1,4(\jjj -1)k+1)&f_0(1,4(\jjj -1)k+2)&0&\dots\\
0&f_1(2,4(\jjj -1)k+2)&f_0(2,4(\jjj -1)k+3)&0&\dots\\
&\ddots&\ddots&\\
&&\ddots&\ddots&\\
&&f_1(4k-2,4\jjj k-2)&f_0(4k-2,4\jjj k-1)&0\\
&&0&f_1(4k-1,4\jjj k-1)&f_0(4k-1,4\jjj k)\\
\endsmallmatrix\right)
$$
and
$$
G^X_\jjj =\left(\smallmatrix g_1(1,4(\jjj -1)k+1)&g_0(1,4(\jjj -1)k+2)&0&\dots\\
0&g_1(2,4(\jjj -1)k+2)&g_0(2,4(\jjj -1)k+3)&0&\dots\\
&\ddots&\ddots&\\
&&\ddots&\ddots&\\
&&g_1(4k-2,4\jjj k-2)&g_0(4k-2,4\jjj k-1)&0\\
&&0&g_1(4k-1,4\jjj k-1)&g_0(4k-1,4\jjj k)\\
\endsmallmatrix\right),
$$
Clearly, we have 
\begin{equation} \label{eq:M^X}
M'''=M^X\big\vert_{X_{1,t}=X_{2,t}=N(k),\,Y_t=1}.
\end{equation}

The evaluation of $\det M^X$ is given in the theorem
below. From the result it is obvious that 
$\det M^X\big\vert_{X_{1,t}=X_{2,t}=N(k),\,Y_t=1}$ is nonzero, and,
thus, also $\det M$.

\begin{theorem} \label{T2}
We have
\begin{multline} \label{eq:CK4}
\det M^X=(-1)^{k-1}4^{2k(4k^2+7k+2)}k^{2k(4k+1)}
\prod _{i=1} ^{4k}(i+1)_{4k-i+1}\\
\times
\prod _{a=1}
^{4k-1}\left(2X_{1,a}-(32k^2+2a-1)Y_a\right)\\
\times
\prod _{1\le a\le b\le 4k-1} ^{}
(2X_{2,b}Y_a-2X_{1,a}Y_b-(2b-2a+1)Y_aY_b),
\end{multline}
where $(\alpha)_k$ is the standard notation for {\it shifted factorials},
$(\alpha)_k:=\alpha(\alpha+1)\cdots(\alpha+k-1)$, $k\ge1$, and $(\alpha)_0:=1$.
\end{theorem}

\begin{remark}\rm
Once having found this theorem,
it is not difficult to prove \eqref{eq:det(A(k))}, by
working out how the coefficients of $P(z)$ resulting from the Ansatz
\eqref{eq:A3}--\eqref{eq:P(t)} are related to the
coefficients of $P(z)$ resulting from the Ansatz \eqref{eq:Ansatz2}. 
Since this is not essential for the
proof of Theorem~\ref{T1}, we leave the details to the reader.
\end{remark}

\medskip
\noindent
{\bf Proof of Theorem~\ref{T2}}.
We follow the ``identification of factors" method as
described in Section~2.4 in \cite{KratBN}.

First we show that $\left(2X_{1,a}-(32k^2+2a-1)Y_a \right)$ divides $\det M^X$,
$a=1,2,\dots,4k-1$. What
has to be proved is that $\det M^X$ vanishes for
$X_{1,a}=(32k^2+2a-1)Y_a/{2}$. 
This can be done by showing that for this choice of $X_{1,a}$ there
is a nontrivial linear combination of the rows of $M^X$. Indeed, if 
$X_{1,a}=(32k^2+2a-1)Y_a/2$ we have
\begin{multline*} \frac {2(X_{2,4k-1}-(N(k)-1)Y_{4k-1})} {(-4)^{k(4k+1)+1}(16k^2+1)\prod
_{\rrrr=1} ^{4k-1}(4\rrrr k+1)}\cdot(\text {row 0 of $M^X$})
\\+
\sum _{\sssss =0} ^{4k}\sum _{\ttttt =0} ^{4k-a-1}\Bigg(\frac {(-1)^{\sssss (k-1)}} {4^{\sssss k}}
\prod _{\rrrr=0} ^{\sssss -1}\frac {4k-1+4\rrrr k} {16k^2+1-4\rrrr k}\kern4cm\\
\cdot
2^\ttttt \prod _{\rrrr=4k-\ttttt } ^{4k-1}\frac {2X_{1,\rrrr}-(32k^2+2\rrrr-1)Y_\rrrr} 
{X_{2,\rrrr-1}-(16k^2+\rrrr-1)Y_{\rrrr-1}}\Bigg)\\
\cdot(\text {row $(16k^2-(4k-1)\sssss -\ttttt -1)$ of $M^X$})=0,
\end{multline*}
as is easy to verify. 

Next we claim that $(2X_{2,b}Y_a-2X_{1,a}Y_b-(2b-2a+1)Y_aY_b)$ 
divides $\det M^X$, $1\le a\le b\le 4k-1$. Let us first impose the
additional restriction that $a<b$. Using the above reasoning,
the claim then follows from the fact that if
$X_{2,b}=\frac {Y_b} {Y_a}X_{1,a}+(2b-2a+1)Y_b/2$ we have
\begin{multline*} 
\sum _{\sssss =0} ^{4k}\sum _{\ttttt =4k-b-1} ^{4k-a-1}
\frac {1} {(-4)^{\sssss k}}\left(\prod _{\rrrr=1} ^{\sssss }\frac 
{2X_{1,a}-(32k^2-8k\rrrr+2a+1)Y_{a}} 
{2X_{1,a}-(64k^2+8k-8k\rrrr+2a+1)Y_{a}}\right)\\
\cdot
4^{b-4k+\ttttt +1}\left(\prod _{\rrrr=4k-\ttttt } ^{b}
\frac {2X_{1,\rrrr}Y_{a}-2X_{1,a}Y_{\rrrr}-(2\rrrr-2a)Y_{a}Y_{\rrrr}} 
{2X_{2,\rrrr-1}Y_{a}-2X_{1,a}Y_{\rrrr-1}-(2\rrrr-2a-1)Y_{a}Y_{\rrrr-1}}\right)\\
\cdot(\text {row $(16k^2-(4k-1)\sssss -\ttttt -1)$ of $M^X$})=0,
\end{multline*}
as is again easy to verify. On the other hand, if $a=b$, then the
same argument shows that 
$(2X_{2,a}-2X_{1,a}-Y_a)$ divides $\det M^X$. It remains to be
checked that also $Y_a$ divides $\det M^X$. Indeed, if $Y_a=0$ then
we have
$$\sum _{\sssss =0} ^{4k}\frac {1} {(-4)^{\sssss k}}
\cdot(\text {row $(16k^2-(4k-1)\sssss -4k+a)$ of $M^X$})=0,
$$
whence $Y_a$ divides $\det M^X$ for $a=1,2,\dots,4k-1$.

These arguments show that the product on the right-hand side of \eqref{eq:CK4}
divides $\det M^X$ as a polynomial in the $X_{1,a}$'s, $X_{2,a}$'s,
and $Y_a$'s.

\medskip
Clearly, the degree in the $X_{1,a}$'s, $X_{2,a}$'s,
and $Y_a$'s of $\det M^X$ is at most $16k^2-1$. But the degree of
the right-hand side of \eqref{eq:CK4} is exactly $16k^2-1$. Therefore we have
proved that 
\begin{multline} \label{eq:CK6}
\det M^X=C_1\prod _{i=1} ^{4k}(i+1)_{4k-i+1}
\prod _{a=1}
^{4k-1}\left(2X_{1,a}-(32k^2+2a-1)Y_a\right)\\
\times
\prod _{1\le a\le b\le 4k-1} ^{}
(2X_{2,b}Y_a-2X_{1,a}Y_b-(2b-2a+1)Y_aY_b),
\end{multline}
where $C_1$ is a constant independent of the $X_{1,a}$'s, $X_{2,a}$'s,
and $Y_a$'s.

In order to determine $C_1$, we compare coefficients of 
\begin{equation} \label{eq:CK7}
X_{1,1}^{4k}X_{1,2}^{4k-1}\cdots X_{1,4k-1}^2Y_1^1Y_2^2\cdots
Y_{4k-1}^{4k-1}
\end{equation}
on both sides of \eqref{eq:CK6}. We claim that the coefficient of this monomial
in $\det M^X$ is equal to $\det M^C$, where $M^C$ is defined exactly
in the same way as $M^X$ (see \eqref{eq:CK3}), except that the definitions of the
functions $f_0,f_1,g_0,g_1$ are replaced by
{\refstepcounter{equation}\label{eq:CK8}}
\alphaeqn
\begin{align} f_0(t,j)&=(N(k)+j)(-4)^k,\label{eq:CK8a}\\
f_1(t,j)&=4(-4)^k,\label{eq:CK8b}\\
g_0(t,j)&=-j,\label{eq:CK8c}\\
g_1(t,j)&=-4.\label{eq:CK8d}
\end{align}
\reseteqn
This is seen as follows. The monomial \eqref{eq:CK7} does not contain any
$X_{2,a}$. Therefore, for finding its coefficient in $\det M^X$, we
may set $X_{2,a}=0$ in $M^X$ for all $a$. 

In which way may the monomial \eqref{eq:CK7} appear in $\det M^X$ (with all $X_{2,a}$
equal to 0)? A typical term in the expansion of $\det M^X$ is the
product of $16k^2$ entries of $M^X$, each from a different row and
column. The monomial \eqref{eq:CK7} contains $X_{1,1}^{4k}$. The variable
$X_{1,1}$ is only found in columns $4\jjj k+1$, $\jjj =0,1,\dots,4k-1$ (and
rows labelled by numbers $\equiv 1$ mod $4k-1$, according to our
labelling scheme). Therefore in a product of entries (each from a
different row and column) which produces a term containing
$X_{1,1}^{4k}$ {\it all\/} the entries from 
columns $4\jjj k+1$ must be ones containing $X_{1,1}$. This explains the above
definitions \eqref{eq:CK8b} and \eqref{eq:CK8d} of $f_1(1,4\jjj k+1)$ and $g_1(1,4\jjj k+1)$,
$\jjj =0,1,\dots,4k-1$, respectively. Moreover, we {\it must\/} generate
the $Y_1$ in \eqref{eq:CK7} from an entry in a column $4\jjj k+2$, for some $\jjj $.
(The variable $Y_1$ is also found in entries in columns $4\jjj k+1$, but
these columns are already taken by our choice of entries which
contain the $X_{1,1}$'s.) This explains the definitions \eqref{eq:CK8a} and \eqref{eq:CK8c}
of $f_0(1,4\jjj k+2)$ and $g_0(1,4\jjj k+2)$, $\jjj =0,1,\dots,4k-1$,
respectively. Next we ask how we can find (in the remaining columns
and rows) entries which contain $X_{1,2}^{4k-1}$. Arguing in an
analogous manner, the variable $X_{1,2}$ only appears in columns
$4\jjj k+2$, $\jjj =0,1,2\dots,4k-1$. One of these columns is already taken
by the entry from which we picked $Y_1$. Therefore in
{\it all\/} the remaining ones we must choose entries 
containing $X_{1,2}$. This explains the 
definitions \eqref{eq:CK8b} and \eqref{eq:CK8d} of $f_1(2,4\jjj k+2)$ and $g_1(2,4\jjj k+2)$,
$\jjj =0,1,\dots,4k-1$, respectively.
Next we consider the term $Y_2^2$ in \eqref{eq:CK7}. It
must come from two entries in columns $4\jjj k+3$, for two different
$\jjj $'s. This explains the definitions \eqref{eq:CK8a} and \eqref{eq:CK8c}
of $f_0(2,4\jjj k+3)$ and $g_0(2,4\jjj k+3)$, $\jjj =0,1,\dots,4k-1$,
respectively. Etc.

The evaluation of $\det M^C$ follows from Lemma~\ref{l7} below
with $X_a=1$ and $Z_a=N(k)$ for $a=1,2,\dots,4k-1$.\QED
\medskip

We consider now a more general determinant than $\det M^C$, the
latter having been
defined through the functions in \eqref{eq:CK8}. Replace these functions by
\begin{align*} f_0(t,j)&=(Z_t+j)(-4)^k,\\
f_1(t,j)&=4(-4)^kX_t,\\
g_0(t,j)&=-j,\\
g_1(t,j)&=-4X_t.
\end{align*}
Let us denote the matrix defined by these functions in the same way
as before $M^C$ by $M^Z$.
Clearly, $M^Z$ specializes to $M^C$ if all $X_t$ are set
equal to 1 and all $Z_t$ to $N(k)$.

The determinant of $M^Z$ evaluates as follows.

\begin{lemma}\label{l7}
We have
\begin{equation} \label{eq:CK9}
\det M^Z=(-1)^{k-1}2^{16k^3+20k^2+14k-1}k^{4k}(4k+1)!
\prod _{a=1} ^{4k-1}\Bigg(X_a^{4k+1-a}
\prod _{b=0} ^{a-1}(Z_a-4bk)\Bigg).
\end{equation}
\end{lemma}
\begin{proof} We proceed in a similar way as in the proof of 
Theorem~\ref{T2}. In the first step we show that the product on the 
right-hand side of \eqref{eq:CK9}
divides $\det M^Z$ as a polynomial in the $X_a$'s and $Z_a$'s. Then,
in the second step,
we compare the degrees of the product and $\det M^Z$. Since
the degree of $\det M^Z$ turns out to be at most the degree of the
product, it then follows that $\det M^Z$ is equal to the
product times some constant which is independent of the $X_t$'s and
$Z_t$'s. Finally, in the third step, this constant is found by
computing the leading coefficient of $\det M^Z$.

\smallskip
{\it Step 1. The product $\prod _{a=1} ^{4k-1}\big(X_a^{4k+1-a}
\prod _{b=0} ^{a-1}(Z_a-4bk)\big)$ divides $\det M^Z$}. 
We start by applying several row and column operations to $\det M^Z$, 
with the final goal
of reducing the size of the determinant. First, for $i=16k^2-1,
16k^2-2,\dots,4k$, in this order, we add $(-4)^k$ times row~$i$ to
row $i-4k+1$. (It should be recalled that, according to our labelling
scheme, we number the rows of $M^Z$ from $0$ to $16k^2-1$.) Thus, we obtain
the determinant of the following matrix:
$$\begin{pmatrix} x_1&x_2&x_3&\hdotsfor2 &x_{4k}\\
F'_1&F'_2&F'_3&\hdotsfor2 &F'_{4k}\\
G_1&F'_1&F'_2&\hdotsfor2 &F'_{4k-1}\\
0&G_2&F'_1&\hdotsfor2 &F'_{4k-2}\\
0&0&G_3&\ddots& &\vdots\\
& &\ddots & \ddots&\ddots&F'_2\\
&& &0 & G_{4k-1}&F'_1\\
& & &\dots & 0&G_{4k}
\end{pmatrix},
$$
where the $x_\ttt $'s, $\ttt =1,2,\dots,4k$, are defined as earlier, and where 
the $(F'_\ttt )$'s and $G_\ttt $'s, $\ttt =1,2,\dots,4k$, are the $(4k-1)\times
(4k)$ matrices 
$$
F'_{\ttt }=\left(\smallmatrix 0&(-4)^{\ttt k}Z_1&0&\dots\\
0&0&(-4)^{\ttt k}Z_2&0&\dots\\
&\ddots&\ddots&\\
&&\ddots&\ddots&\\
&&&0&(-4)^{\ttt k}Z_{4k-2}&0\\
&&&0&0&(-4)^{\ttt k}Z_{4k-1}\\
\endsmallmatrix\right)
$$
and
$$
G_\ttt =\left(\smallmatrix -4X_1&-4(\ttt -1)k-2&0&\dots\\
0&-4X_2&-4(\ttt -1)k-3&0&\dots\\
&\ddots&\ddots&\\
&&\ddots&\ddots&\\
&&0&-4X_{4k-2}&-4\ttt k-1&0\\
&&&0&-4X_{4k-1}&-4\ttt k\\
\endsmallmatrix\right).
$$

Next we ``make" the submatrices $G_\ttt $, $\ttt =1,2,\dots,4k$, to diagonal
matrices, by subtracting $(j+1)/4X_j$ times column~$j$ from column
$j+1$, $j=1,2,\dots,4k-1$, $(j+1)/4X_{j-4k}$ times column~$j$ from
column $j+1$, $j=4k+1,4k+2,\dots,8k-1$, 
\dots, and $(j+1)/4X_{j-16k^2+4k}$ times column~$j$ from
column $j+1$, $j=16k^2-4k+1,16k^2-4k+2,\dots,16k^2-1$. After these
operations we obtain the determinant of the matrix
$$\begin{pmatrix} x_1&x_2&x_3&\hdotsfor2 &x_{4k}\\
F_{1,1}&F_{2,2}&F_{3,3}&\hdotsfor2 &F_{4k,4k}\\
G&F_{1,2}&F_{2,3}&\hdotsfor2 &F_{4k-1,4k}\\
0&G&F_{1,3}&\hdotsfor2 &F_{4k-2,4k}\\
0&0&G&\ddots&&\vdots\\
& &\ddots & \ddots&\ddots&F_{2,4k}\\
&& &0 & G&F_{1,4k}\\
& & &\dots & 0&G
\end{pmatrix},
$$
where $G$ and the $F_{\ssss ,\tttt }$'s, 
$1\le \ssss \le \tttt \le 4k$, are the $(4k-1)\times
(4k)$ matrices 
$$
G=\left(\smallmatrix -4X_1&0&0&\dots\\
0&-4X_2&0&0&\dots\\
&\ddots&\ddots&\ddots\\
&&\ddots&\ddots&\ddots&\\
&&&0&-4X_{4k-2}&0&0\\
&&&&0&-4X_{4k-1}&0\\
\endsmallmatrix\right)
$$
and
$$
F_{\ssss ,\tttt }=\left(\smallmatrix
0&f^{(\ssss ,\tttt )}_{1,2}&f^{(\ssss ,\tttt )}_{1,3}&\dots&&f^{(\ssss ,\tttt )}_{1,4k}\\
0&0&f^{(\ssss ,\tttt )}_{2,3}&f^{(\ssss ,\tttt )}_{2,4}&\dots&f^{(\ssss ,\tttt )}_{2,4k}\\
&\ddots&\ddots&&&\vdots\\
&&\ddots&\ddots&&\vdots\\
&&&0&f^{(\ssss ,\tttt )}_{4k-2,4k-1}&f^{(\ssss ,\tttt )}_{4k-2,4k}\\
&&&0&0&f^{(\ssss ,\tttt )}_{4k-1,4k}\\
\endsmallmatrix\right),
$$
where 
$$f^{(\ssss ,\tttt )}_{\iiii \jjjj }=\frac {(-4)^{\ssss k-\jjjj +\iiii +1}Z_\iiii \,(4(\tttt -1)k+\iiii +2)_{\jjjj -\iiii -1}}
{X_{\iiii +1}X_{\iiii +2}\cdots X_{\jjjj -1}}.
$$

Now we eliminate the last columns in $F_{\ssss \tttt }$ for $1\le \ssss <\tttt \le 4k$.
We start by eliminating the last column of $F_{1,4k}$. We do this by
adding $f^{(1,4k)}_{\iiii ,4k}/4X_\iiii $ times column $16k^2-8k+\iiii $ to column
$16k^2$, $\iiii =1,2,\dots,4k-1$. This makes all the entries in the last
column which are in rows $16k^2-8k+2,\dots,16k^2-4k-1,16k^2-4k$ zero,
whereas the entries in the last column in rows
$16k^2-12k+3,\dots,16k^2-8k+1$ are modified. Next we eliminate these
entries in a similar fashion, by using the columns $16k^2-12k+\iiii $,
$\iiii =1,2,\dots,4k-1$, etc. In the end all the entries in the last
column in rows $4k,4k+1,\dots,16k^2-1$ will be zero, whereas the
entries in the last column in rows $1,2,\dots,4k-1$ will have been
(significantly) modified. An analogous procedure is applied to
eliminate the entries in the last columns of $F_{\ssss ,4k-1}$,
$\ssss =1,2,\dots,4k-2$. Just to mention the first step: We
add $f^{(1,4k-1)}_{\iiii ,4k}/4X_\iiii $ times column $16k^2-12k+\iiii $ to column
$16k^2-4k$, $\iiii =1,2,\dots,4k-1$. This makes all the entries in 
column $16k^2-4k$ which are in rows 
$16k^2-12k+3,\dots,16k^2-8k,16k^2-8k+1$ zero,
whereas the entries in the last column in rows
$16k^2-16k+4,\dots,16k^2-12k+2$ are modified. Etc.

The advantage after having done all this is that now all the entries
in columns $4\jjj k$, $\jjj =1,2,\dots,4k$, are zero except for entries in
rows $0,1,\dots,4k-1$.
This fact, and the fact that the submatrices $G$ are diagonal
matrices (of rectangular form) with last column consisting entirely
of zeroes, makes it possible
to reduce the determinant of the (new) matrix significantly. For
$i=16k^2-1,16k^2-2,\dots,4k$ we may
expand the determinant with respect to row~$i$, in this order. If the
details are worked out, then we see that our original determinant
$\det M^Z$ is equal to
$$(-4)^{16k^2-4k}\prod _{a=1} ^{4k-1}X_a^{4k}
\det\begin{pmatrix} 
u_1&u_2&\dots&u_{4k}\\
e_{1,1}&e_{1,2}&\dots&e_{1,4k}\\
\hdotsfor4\\
e_{4k-1,1}&e_{4k-1,2}&\dots&e_{4k-1,4k}\\
\end{pmatrix},
$$
where $u_\ell$ is given by \eqref{eq:CK2}, and where $e_{ij}$ is given by
\begin{multline*}
\sum _{\rrr=1} ^{j}\underset{ j+1=n_0> n_1>n_2>\dots>n_\rrr=1}
{\sum _{4k=t_0> t_1>\dots>t_\rrr=i} ^{}}
(-1)^{\rrr-1}(-4)^{jk-4k+i+1}Z_{t_1}Z_{t_2}\cdots Z_{t_\rrr}\\
\times \frac {\prod _{\nu=0} ^{\rrr-1}(4k(n_\nu-1)-4k+t_{\nu+1}+2)
_{t_\nu-t_{\nu+1}-1}} {X_{i+1}X_{i+2}\cdots X_{4k-1}}.
\end{multline*}
Clearly, we may extract $(-4)^{jk}$ from column~$j$,
$j=1,2,\dots,4k$, and\linebreak $(-4)^{-4k+i+1}/{X_{i+1}X_{i+2}\cdots X_{4k-1}}$ from
row~$i$, $i=1,2,\dots,4k-1$ (still using our nonstandard labelling
scheme where the rows are numbered $0,1,\dots,4k-1$), so that we
obtain the expression
\begin{multline} \label{eq:CK10}
(-4)^{16k^2-4k+2k^2(4k+1)-(2k-1)(4k-1)}\prod _{a=1}
^{4k-1}X_a^{4k+1-a}\\
\times
\det\begin{pmatrix} 
\tilde u_1&\tilde u_2&\dots&\tilde u_{4k}\\
\tilde e_{1,1}&\tilde e_{1,2}&\dots&\tilde e_{1,4k}\\
\hdotsfor4\\
\tilde e_{4k-1,1}&\tilde e_{4k-1,2}&\dots&\tilde e_{4k-1,4k},\\
\end{pmatrix},
\end{multline}
where $\tilde u_\ell$ is given by 
$$(-1)^{\ell-1}(-4)^{k}8k(4k+1)\left(\prod _{i=1} ^{4k-\ell}(4ik-1)\right)
\left(\prod _{i=1} ^{\ell-1}(4ik+1)\right),$$
and where $\tilde e_{ij}$ is given by
\begin{multline} \label{eq:CK11}
\sum _{\rrr=1} ^{j}\underset {j+1=\tilde n_0> \tilde n_1>\tilde n_2>\dots>\tilde
n_\rrr=1}
{\sum _{4k=t_0> t_1>\dots>t_\rrr=i} ^{}}
(-1)^{\rrr-1}Z_{t_1}Z_{t_2}\cdots Z_{t_\rrr}\\
\times  {\prod _{\nu=0} ^{\rrr-1}(4k(\tilde n_\nu-1)-4k+t_{\nu+1}+2)
_{t_\nu-t_{\nu+1}-1}}.
\end{multline}
{}From \eqref{eq:CK10} it is abundantly clear that $\prod _{a=1} ^{4k-1}X_a^{4k+1-a}$
divides $\det M^Z$ (which, after all, is equal to \eqref{eq:CK10}). It remains to
show that also $\prod _{a=1} ^{4k-1}\prod _{b=0} ^{a-1}(Z_a-4bk)$
divides $\det M^Z$.

Let $a$ and $b$, $1\le a\le 4k-1$, $0\le b\le a-1$, be given. We want
to show that $Z_a-4bk$ divides $\det M^Z$. We will show the
equivalent fact that the rows of $\det M^Z\big\vert_{Z_a=4bk}$ are
linearly dependent. The crucial observation, from which this claim
follows easily, is that the entries $\tilde e_{ij}$,
$j=1,2,\dots,4k$, in row~$i$ of the determinant in \eqref{eq:CK10} are given by a
polynomial in $j$, $p_i(j)$ say, of degree $(4k-i-1)$ and with leading
coefficient (i.e., coefficient of $j^{4k-i-1}$) equal to
\begin{equation} \label{eq:CK12}
\sum _{\rrr=1} ^{j}{\sum _{4k=t_0> t_1>\dots>t_\rrr=i} ^{}}
(-1)^{\rrr-1}Z_{t_1}Z_{t_2}\cdots Z_{t_\rrr}\frac {(4k)^{4k-i-\rrr}}
{\prod _{\nu=0} ^{\rrr-1}(t_\nu-i)}.
\end{equation}
This is seen as follows. The summand in \eqref{eq:CK11} is a polynomial in 
$j,\tilde n_1,\tilde n_2,\dots,\tilde n_{\rrr-1}$
of multidegree
$(4k-t_1-1,t_1-t_2-1,t_2-t_3-1,\dots,t_{\rrr-1}-t_\rrr-1)$ (i.e., the
degree in $j$ is $4k-t_1-1$, the degree in $\tilde n_1$ is
$t_1-t_2-1$, etc.). Because of the fact that (for fixed $u$ and varying $v$)
$\sum _{\gamma=u} ^{v}\gamma^e$ is a polynomial in $v$ of degree $e+1$ with leading
coefficient (i.e., coefficient of $v^{e+1}$) equal to $1/(e+1)$,
successive summation over $\tilde n_{\rrr-1}$, $\tilde n_{\rrr-2}$, \dots,
$\tilde n_1$
yields the claimed facts.

Under the specialization $Z_a=4bk$, it is seen by ``inspection" that
\eqref{eq:CK12} with $i$ replaced by $a-b$ vanishes, because
the summand corresponding to $t_0>\dots>t_{\omega-1}>t_\omega=a>t_{\omega+1}
>\dots>t_r$
cancels with the summand corresponding to $t_0>\dots>t_{\omega-1}>
t_{\omega+1}>\dots>t_r$. Hence, the polynomial 
$p_{a-b}(j)$ has degree (at most) $4k-a+b-2$ (instead of $4k-a+b-1$).
Consequently, if $Z_a=4bk$ then 
the entries in rows $a-b$, $a-b+1$, \dots, $4k-1$ are given by
polynomials in $j$ (to wit: $j$ denoting the column index of the
entries) of respective degrees $4k-a+b-2$ (!), $4k-a+b-2$,
$4k-a+b-3$, \dots, $2$, $1$. These are $4k-a+b$ polynomials, all of degree
at most $4k-a+b-2$. It follows that there must be a nontrivial linear
combination of these polynomials that vanishes. Hence, the rows 
$a-b$, $a-b+1$, \dots, $4k-1$ are linearly dependent, which, in turn,
implies that the determinant in \eqref{eq:CK10} (and, thus, also $\det M^Z$)
vanishes for $Z_a=4bk$, $1\le a\le 4k-1$, $0\le b\le a-1$.

\smallskip
{\it Step 2. Comparison of degrees}.
Clearly, the degree of $\det M^Z$ as a polynomial in the $X_t$'s and
$Z_t$'s is at most $16k^2-1$, whereas the degree of the product on
the right-hand side of \eqref{eq:CK9} is exactly $16k^2-1$. Hence, we
have
$$\det M^Z=C_2\prod _{a=1} ^{4k-1}\big(X_a^{4k+1-a}
\prod _{b=0} ^{a-1}(Z_a-4bk)\big),$$
where $C_2$ is a constant independent of the $X_t$'s and $Z_t$'s.

\smallskip
{\it Step 3. Computation of the leading coefficient}.
In order to determine $C_2$, we determine the coefficient of 
$\prod _{a=1} ^{4k-1}X_a^{4k+1-a}Z_a^a$ in the expansion of $\det
M^Z$. By arguments similar to those at the end of the proof of
Theorem~\ref{T2}, it is seen that this coefficient is given by the determinant
of the following matrix, which we denote by $M^L$. It is defined
exactly in the same way as $M^X$ (see \eqref{eq:CK3}), except that the definitions of the
functions $f_0,f_1,g_0,g_1$ are replaced by
\begin{align*} f_0(t,j)&=(-4)^k,\\
f_1(t,j)&=4(-4)^k,\\
g_0(t,j)&=0,\\
g_1(t,j)&=-4.
\end{align*}
By expanding this determinant with respect to row~0, we obtain
\begin{equation} \label{eq:CK13}
\det M^L=\sum _{\jj =1} ^{4k}(-1)^{\jj -1}u_\jj \det M^L_\jj ,
\end{equation}
where $u_\jj $ is, as earlier, given by \eqref{eq:CK2}, and $M^L_\jj $ is the
matrix arising from $M^L$ by deleting row~0 and column $4\jj k$.

Let $\jj $, $1\le \jj \le 4k$, be fixed. We will next compute $\det M^L_\jj $.
When we built $M^L_\jj $ from $M^L$, we deleted in particular row~0. 
Therefore we will now switch
to the usual labelling scheme for rows and columns of a matrix, i.e.,
we will subsequently not only label the columns by $1,2,\dots$ but
also the rows.

If $\jj <4k$, then
we expand $\det M^L_\jj $ with respect to the last $4k-1$ rows and then
with respect to the last column. Since these rows and this column
contain only one nonzero entry, we obtain some multiple of the
determinant of a
$(16k^2-4k-1)\times (16k^2-4k-1)$ matrix. If $\jj <4k-1$, then we
continue by expanding the (now reduced) determinant with respect to
the last $4k-2$ rows and then with respect to the last 2 columns. We
continue in the same manner until we have reduced $\det M^L_\jj $ to the
determinant of a
$(4\jj k-1)\times (4\jj k-1)$ matrix, more precisely, until we arrive
at
$$\det M^L_\jj =4^{\binom {4k}2-\binom \jj 2}(-4)^{k\binom {4k-\jj +1}2}
\det M'_\jj ,$$
where $M'_\jj $ is the matrix
$$\begin{pmatrix} 
F&0&0&\hdotsfor2 &0\\
G&F&0&\hdotsfor2 &0\\
0&G&F&\hdotsfor2 &0\\
\vdots & & \ddots&\ddots&&\vdots\\
0&\dots  &0 & G&F&0\\
0&\hdotsfor2   &0 & G&F'\\
0&\hdotsfor3   & 0&U\\
\end{pmatrix},\quad \quad \text {($\jj $ occurrences of $F$)},
$$
with $F$ and $G$ the $(4k-1)\times
(4k)$ matrices 
$$
F=\begin{pmatrix} 4(-4)^k&(-4)^k&0&\dots\\
0&4(-4)^k&(-4)^k&0&\dots\\
&\ddots&\ddots&\ddots\\
&&\ddots&\ddots&\ddots&\\
&&&0&4(-4)^k&(-4)^k&0\\
&&&&0&4(-4)^k&(-4)^k\\
\end{pmatrix}
$$
and
$$
G=\begin{pmatrix} -4&\hphantom{-}0&0&\dots\\
\hphantom{-}0&-4&0&0&\dots\\
&\ddots&\ddots&\ddots\\
&&\ddots&\ddots&\ddots&\\
&&&\hphantom{-}0&-4&\hphantom{-}0&0\\
&&&&\hphantom{-}0&-4&0\\
\end{pmatrix},
$$
$F'$ the $(4k-1)\times(4k-1)$ matrix which arises from $F$ by
deleting its last column,
and $U$ the $(\jj -1)\times (4k-1)$ matrix
$$U=\begin{pmatrix} 
-4&\hphantom{-}0&0&\hdotsfor3 &0\\
\hphantom{-}0&-4&0&\hdotsfor3 &0\\
\hphantom{-}\vdots&&\ddots&\ddots& &&\vdots\\
\hphantom{-}0&\hdotsfor2 &-4&0&\dots &0
\end{pmatrix}.
$$
We continue by expanding $M'_\jj $ with respect to the last $\jj -1$ rows.
Thus we obtain
\begin{equation} \label{eq:CK14}
\det M^L_\jj =4^{\binom {4k}2-\binom \jj 2+\jj -1}(-4)^{k\binom {4k-\jj +1}2}
\det \overline M_\jj 
\end{equation}
for $\det M^L_\jj $,
where $\overline M_\jj $ is the matrix
$$\begin{pmatrix} 
F&0&0&\hdotsfor2 &0\\
G&F&0&\hdotsfor2 &0\\
0&G&F&\hdotsfor2 &0\\
\vdots & & \ddots&\ddots&&\vdots\\
0&\dots &0 & G&F&0\\
0&\hdotsfor2  &0 & G&V\\
\end{pmatrix},\quad \quad \text {($\jj -1$ occurrences of $F$)},
$$
with $V$ the $(4k-1)\times (4k-\jj )$ matrix
\begin{equation} \label{eq:CK15}
V=\left(\begin{matrix} 0&0&\hdotsfor3 &0\\
\hdotsfor6 \quad  \\
0&0&\hdotsfor3 &0\\
(-4)^k&0&\dots\\
4(-4)^k&(-4)^k&0&\dots\\
0&4(-4)^k&(-4)^k&0&\dots\\
&\ddots&\ddots&\ddots\\
&&\ddots&\ddots&\ddots&\\
&&&0&4(-4)^k&(-4)^k\\
&&&&0&4(-4)^k\\
\end{matrix}\right).
\end{equation}

Next we prepare for a reduction from the top of $\overline M_\jj $. We
subtract 4 times column $j$ from column $j-1$, $j=4(\jj -1)k,
4(\jj -1)k-1,\dots,4k(\jj -2)+2$, $j=4k(\jj -2),
4k(\jj -2)-1,\dots,4k(\jj -3)+2$, \dots, $j=4k,
4k-1,\dots,2$, in this order. Thus $\det M'_\jj $ is converted to
$\det M''_\jj $, where $M''_\jj $ is the matrix
$$\begin{pmatrix} 
F''&0&0&\hdotsfor3 &0\\
G''&F''&0&\hdotsfor3 &0\\
0&G''&F''&\hdotsfor3 &0\\
0&0&G''&\ddots& &&\vdots\\
\vdots& & \ddots& \ddots&\ddots&&0\\
0&\hdotsfor2  &0& G''&F''&0\\
0&\hdotsfor3   &0& G''&V
\end{pmatrix},
$$
with $F''$ and $G''$ the $(4k-1)\times
(4k)$ matrices 
$$
F''=\begin{pmatrix} 0&(-4)^k&0&\dots\\
0&0&(-4)^k&0&\dots\\
&\ddots&\ddots&\ddots\\
&&\ddots&\ddots&\ddots&\\
&&&0&0&(-4)^k&0\\
&&&&0&0&(-4)^k\\
\end{pmatrix}
$$
and
$$
G''=\begin{pmatrix} -4&\hphantom{-}0&0&\dots\\
(-4)^2&-4&0&0&\dots\\
\vdots&\ddots&\ddots&\ddots\\
\vdots&&\ddots&\ddots&\ddots&\\
(-4)^{4k-2}&\hdotsfor2 &(-4)^2&-4&\hphantom{-}0&0\\
(-4)^{4k-1}&\hdotsfor2 &(-4)^3&(-4)^2&-4&0\\
\end{pmatrix}.
$$

Our next goal is to ``push" the (nonzero) entries in columns $4tk+1$,
$t=0,1,\dots,\jj -2$, down to rows
$(4k-1)(\jj -1)+1,\dots,(4k-1)\jj -1,(4k-1)\jj $. (This is similar to what we
did in Step~1 when we ``pushed" all the nonzero entries in columns 
$4tk$, $t=1,2,\dots,4k$ up to rows $0,1,\dots,4k-1$.)
In order to achieve this for the 1-st column, we add
$$\sum _{r=1} ^{\jj -2}\sum _{s=2} ^{4k}
(-1)^r(-4)^{-r k+s-1}\binom {s-1}{r-1}
\cdot(\text {column $(4kr+s)$})$$
to column 1. Similarly, in order to achieve this for the $(4k+1)$-st 
column, we add
$$\sum _{r=2} ^{\jj -2}\sum _{s=2} ^{4k}
(-1)^{r-1}(-4)^{-(r-1) k+s-1}\binom {s-1}{r-2}
\cdot(\text {column $(4kr+s)$})$$
to column $4k+1$. Etc. As a result, the determinant $\det M''_\jj $ is
converted to the determinant of the matrix
\begin{equation} \label{eq:CK16}
\begin{pmatrix} 
F''&0&0&\hdotsfor3 &0\\
G'&F''&0&\hdotsfor3 &0\\
0&G'&F''&\hdotsfor3 &0\\
0&0&G'&\ddots& &&\vdots\\
\vdots& & \ddots& \ddots&\ddots&&0\\
0&\hdotsfor2  &0& G'&F''&0\\
H_1&H_2&\hdotsfor2  &H_{\jj -2}&G''&V
\end{pmatrix},
\end{equation}
where $F''$ and $G''$ are as before, $G'$ is the $(4k-1)\times (4k)$
matrix
$$
G'=\begin{pmatrix} 0&\hphantom{-}0&0&\dots\\
0&-4&0&0&\dots\\
\vdots&\ddots&\ddots&\ddots\\
\vdots&&\ddots&\ddots&\ddots&\\
0&(-4)^{4k-3}&\dots&(-4)^2&-4&\hphantom{-}0&0\\
0&(-4)^{4k-2}&\dots&(-4)^3&(-4)^2&-4&0\\
\end{pmatrix},
$$
and $H_t$, $t=1,2,\dots,\jj -2$, is a $(4k-1)\times (4k)$ matrix
with all entries equal to 0, except for the entries in the first
column. To be precise, the entry in the first column and row $s$ of
$H_t$ is given by
\begin{equation} \label{eq:CK17}
(-1)^{\jj -1-t}(-4)^{-(\jj -1-t)k+s}\binom {s-1}{\jj -1-t}.
\end{equation}
It should be noted that the entries in the first column of $G''$
($G''$ appearing at the bottom of the matrix \eqref{eq:CK16}, as do the matrices
$H_t$) is given by \eqref{eq:CK17} with $t=\jj -1$.

Now everything is prepared for the reduction. We expand the
determinant of \eqref{eq:CK16} with respect to rows $1,2,\dots,(4k-1)(\jj -1)$.
This reduces the determinant of \eqref{eq:CK16} to
$$(-1)^{\binom \jj 2}(-4)^{(\jj -1)(4k-1)k}\det \tilde M_\jj ,$$
where $\tilde M_\jj $ is a $(4k-1)\times(4k-1)$ matrix of the form
$$\begin{pmatrix} N&V\end{pmatrix},$$
with the $(s,t)$-entry of $N$ being given by \eqref{eq:CK17},
$s=1,2,\dots,4k-1$, $t=1,2,\dots,\jj -1$, and $V$ the $(4k-1)\times(4k-\jj )$
matrix from above. If we substitute all this in \eqref{eq:CK14}, we obtain that
$\det M^L_\jj $ is equal to
\begin{equation} \label{eq:CK18}
(-1)^{\binom \jj 2}
4^{\binom {4k}2-\binom \jj 2+\jj -1}(-4)^{k\binom {4k-\jj +1}2+(\jj -1)(4k-1)k}
\det \tilde M_\jj .
\end{equation}
The submatrix $V$ of $\tilde M_\jj $ is almost diagonal. Subtraction of 4
times row $s$ from row $s+1$ in $\tilde M_\jj $, $s=\jj -1,\jj ,\dots,4k-2$, will
transform it into a completely diagonal matrix (namely into the
matrix on the right-hand side of \eqref{eq:CK15} with all entries $4(-4)^k$
replaced by 0). As a side effect, this will turn the
$(4k-1,1)$-entry of $\tilde M_\jj $ into 
$$(-1)^{\jj -2}(-4)^{-(\jj -2)k+4k-1}\binom {4k-1}{\jj -1}.$$
As is easily seen, the determinant of the in this way modified
matrix, $M^*_\jj $ say, is
\begin{multline*}
(-1)^{4k-\jj +\binom {\jj -1}2}(-4)^{(4k-\jj )k}\cdot
\left(\text {$(4k-1,1)$-entry of $M^*_\jj $}\right)\\
\times
\prod _{s=1} ^{\jj -2}\left(\text {$(s,\jj -s)$-entry of $M^*_\jj $}\right),
\end{multline*}
or, explicitly,
$$(-1)^\jj (-4)^{(4k-\jj )k-(\jj -2)k+4k-1-\binom {\jj -2}2 k+\binom {\jj -1}2}\binom {4k-1}{\jj -1}.
$$
Substitution of the above in \eqref{eq:CK18} yields that
the determinant $\det M^L_\jj $ is equal to
$$(-1)^{\jj k}
4^{ 8k^3+ 10k^2 - \jj k + 2k-1}
\binom {4k-1}{\jj -1}.
$$

Now we substitute this in \eqref{eq:CK13}. We obtain that $\det M^L$ is equal to
\begin{multline*}
\sum _{\jj =1} ^{4k}
(-1)^{k} 4^{ 8k^3+ 10k^2  + 3k}
\binom {4k-1}{\jj -1}
2k(4k+1)\left(\prod _{i=1} ^{4k-\jj }(4ik-1)\right)
\left(\prod _{i=1} ^{\jj -1}(4ik+1)\right)\\
=(-1)^k 2^{16k^3+20k^2+14k-1}k^{4k}(4k+1)(4k-1)!\,
\sum _{\jj =1} ^{4k}\binom {\frac {1} {4k}-1}{4k-\jj }
\binom {-\frac {1} {4k}-1}{\jj -1}.
\end{multline*}
The sum is readily evaluated by means of the Chu--Vandermonde
summation (see e.g.\ \cite[Sec.~5.1, (5.27)]{GrKPAA}), so that we
obtain
\begin{multline*}
(-1)^k 2^{16k^3+20k^2+14k-1}k^{4k}(4k+1)(4k-1)!\,
\binom {-2}{4k-1}\\
=
(-1)^{k-1} 2^{16k^3+20k^2+14k-1}k^{4k}(4k+1)!\,.
\end{multline*}
Since the coefficient of $\prod _{a=1} ^{4k-1}X_a^{4k+1-a}Z_a^a$ in
the expression on the right-hand side of \eqref{eq:CK9} is exactly the same, we
have completed the proof of the lemma.\QED
\end{proof}


\bigskip


\begin{thebibliography}{1}

\bigskip
\small

\bibitem{Bellard}
Fabrice Bellard, {\em $\pi$ page},
{\tt http://fabrice.bellard.free.fr/pi/}.

\bibitem{BaBoPl}
D. Bailey, P. Borwein, S. Plouffe, {\em On the rapid computation of 
various polylogarithmic constants}, Math.\ Comp.\ {\bf 66} (1997), 
903--913. 

\bibitem{Gosper}
R. Wm. Gosper, unpublished research announcement, 1974.

\bibitem{GrKPAA}
R. L. Graham, D. E. Knuth and O. Patashnik,
{\it Concrete Mathematics}, Addison-Wesley, Reading, Massachusetts,
1989.

\bibitem{KratBN}
C.    Krattenthaler, {\it Advanced determinant calculus},
S\'eminaire Lotharingien Combin.\ {\bf 42} (``The Andrews Festschrift")
(1999), Article~B42q, 67~pp.

\bibitem{Plouffe}
S. Plouffe, {\em On the computation of the $n$'th decimal digit
of various transcendental numbers}, manuscript, 1996, available at {\tt 
http://www.lacim.uqam.ca/plouffe/}.

\end{thebibliography}
\end{document}